\documentclass[reqno]{amsart}
\usepackage{setspace,amssymb}
\usepackage{ifpdf}
\ifpdf
 \usepackage[hyperindex]{hyperref}
\else
 \expandafter\ifx\csname dvipdfm\endcsname\relax
 \usepackage[hypertex,hyperindex]{hyperref}
 \else
 \usepackage[dvipdfm,hyperindex]{hyperref}
 \fi
\fi
\allowdisplaybreaks[4]
\numberwithin{equation}{section}
\theoremstyle{plain}
\newtheorem{thm}{Theorem}[section]

\theoremstyle{remark}
\newtheorem{rem}{Remark}[section]
\newcommand{\bell}{\textup{B}}
\DeclareMathOperator{\td}{d}

\begin{document}

\title[Diagonal recurrence relations for Stirling numbers]
{Diagonal recurrence relations for the Stirling numbers of the first kind}

\author[F. Qi]{Feng Qi}

\address[F. Qi]{Department of Mathematics, College of Science, Tianjin Polytechnic University, Tianjin City, 300387, China}
\email{\href{mailto: F. Qi <qifeng618@gmail.com>}{qifeng618@gmail.com}, \href{mailto: F. Qi <qifeng618@hotmail.com>}{qifeng618@hotmail.com}, \href{mailto: F. Qi <qifeng618@qq.com>}{qifeng618@qq.com}}
\urladdr{\url{http://qifeng618.wordpress.com}}

\begin{abstract}
In the paper, the author presents diagonal recurrence relations for the Stirling numbers of the first kind. As by-products, the author also recovers three explicit formulas for special values of the Bell polynomials of the second kind.
\end{abstract}



\thanks{This paper was typeset using \AmS-\LaTeX}

\maketitle

\section{Introduction}

In combinatorics, the Bell polynomials of the second kind, or say, the partial Bell polynomials, denoted by $\bell_{n,k}(x_1,x_2,\dotsc,x_{n-k+1})$ for $n\ge k\ge0$, are defined by
\begin{equation}
\bell_{n,k}(x_1,x_2,\dotsc,x_{n-k+1})=\sum_{\substack{1\le i\le n,\ell_i\in\{0\}\cup\mathbb{N}\\ \sum_{i=1}^ni\ell_i=n\\
\sum_{i=1}^n\ell_i=k}}\frac{n!}{\prod_{i=1}^{n-k+1}\ell_i!} \prod_{i=1}^{n-k+1}\Bigl(\frac{x_i}{i!}\Bigr)^{\ell_i}.
\end{equation}
See~\cite[p.~134, Theorem~A]{Comtet-Combinatorics-74}. For more information on the Bell polynomials in general and Dyck paths in particular, please look at the papers~\cite{Mansour-Sun-DAM-2008} and~\cite{Mansour-Sun-AJC-2008} and plenty of references therein.
\par
In mathematics, the Stirling numbers arise in a variety of combinatorics problems and were introduced by James Stirling in the eighteen century. There are two different kinds of the Stirling numbers. The Stirling numbers of the first kind $s(n,k)$, which are also called the signed the Stirling numbers of the first kind, may be generated by
\begin{equation}
\frac{[\ln(1+x)]^k}{k!}=\sum_{n=k}^\infty s(n,k)\frac{x^n}{n!},\quad |x|<1. \label{gen-funct-3}
\end{equation}
The mathematical meaning of the unsigned Stirling numbers of the first kind $(-1)^{n-k}s(n,k)$ can be interpreted as the number of permutations of $\{1,2,\dotsc,n\}$ with $k$ cycles.
\par
Several ``triangular'', ``horizontal'', and ``vertical'' recurrence relations for the Stirling numbers of the first kind $s(n,k)$ were listed in~\cite[pp.~214--215, Theorems~A,~B, and~C]{Comtet-Combinatorics-74} as
\begin{align}
s(n,k)&=s(n-1,k-1)-(n-1)s(n-1,k), \label{triangular-rel}\\
(n-k)s(n,k)&=\sum_{k+1\le\ell\le n}(-1)^{\ell-k}\binom{\ell}{k-1}s(n,\ell), \label{horiz-rel-1}\\
s(n,k)&=\sum_{k\le\ell\le n}s(n+1,\ell+1)n^{\ell-k}, \label{horiz-rel-2}\\
ks(n,k)&=\sum_{k-1\le\ell\le n-1}(-1)^{n-\ell-1}\binom{n}{\ell}s(\ell,k-1), \label{vert-rel-1}\\
s(n+1,k+1)&=\sum_{k\le\ell\le n}(-1)^{\ell-1}\prod_{q=1}^{n-\ell}(\ell+q)s(\ell,k), \label{vert-rel-2}
\end{align}
where as usual the empty product means $1$.
\par
The aim of this paper is to present, basing on an integral representation for the Stirling numbers of the first kind $s(n,k)$, making use of Fa\`a di Bruno formula, and utilizing properties of the Bell polynomials of the second kind $\bell_{n,k}$, diagonal recurrence relations for the Stirling numbers of the first kind $s(n,k)$. As by-products, three explicit formulas for special values of the Bell polynomials of the second kind $\bell_{n,k}$ are recovered.
\par
The main results may be formulated in the following theorem.

\begin{thm}\label{s(n-k)=Bell(n-k)-thm}
For $n\ge k\ge1$, we have
\begin{gather}\label{Bell-Stir1st=eq}
\bell_{n,k}\biggl(\frac{1!}2,\frac{2!}3,\dotsc, \frac{(n-k+1)!}{n-k+2}\biggr)
=(-1)^{n-k}\frac1{k!}\sum_{m=1}^k(-1)^m\frac{\binom{k}{m}}{\binom{n+m}{n}}s(n+m,m),\\
\bell_{n,k}(0,1!,\dotsc,(n-k)!)
=(-1)^{n-k}\binom{n}{k} \sum_{m=0}^k(-1)^{m}\frac{\binom{k}{m}}{\binom{n-m}{n-k}}s(n-m,k-m), \label{Bell-Stir1n!=eq}
\end{gather}
and
\begin{align}\label{s(n-k)=s(n-k)-id}
s(n,k)&=(-1)^{k}\sum_{m=1}^{n}(-1)^{m}\sum_{\ell=k-m}^{k-1}(-1)^{\ell} \binom{n}{\ell}\binom{\ell}{k-m} s(n-\ell,k-\ell)\\ \label{1stirling-diagonal-eq}
&=(-1)^{n-k}\sum_{\ell=0}^{k-1}(-1)^\ell\binom{n}{\ell} \binom{\ell-1}{k-n-1}s(n-\ell,k-\ell),
\end{align}
where the conventions that
\begin{equation*}
\binom00=1,\quad \binom{-1}{-1}=1, \quad \text{and}\quad \binom{p}q=0
\end{equation*}
for $p\ge0>q$ are adopted in~\eqref{1stirling-diagonal-eq}.
\end{thm}

\section{Proof of Theorem~\ref{s(n-k)=Bell(n-k)-thm}}

Recently, three integral representations for the Stirling numbers of the first kind $(-1)^{n-k}s(n,k)$ were discovered in~\cite{1st-Sirling-Number-2012.tex}. The first one among them, \cite[Theorem~2.1]{1st-Sirling-Number-2012.tex}, reads that, for $1\le k\le n$,
\begin{equation}\label{n-times-diriv}
s(n,k)=\binom{n}{k}\lim_{x\to0}\frac{\td^{n-k}}{\td x^{n-k}} \biggl\{\biggl[\int_0^\infty\biggl(\int_{1/e}^1 t^{xu-1}\td t\biggr)e^{-u}\td u\biggr]^k\biggr\}.
\end{equation}
\par
In combinatorial analysis, Fa\`a di Bruno formula plays an important role and may be described in terms of the Bell polynomials of the second kind $\bell_{n,k}$ by
\begin{equation}\label{Bruno-Bell-Polynomial}
\frac{\td^n}{\td t^n}f\circ h(t)=\sum_{k=1}^nf^{(k)}(h(t)) \bell_{n,k}\bigl(h'(t),h''(t),\dotsc,h^{(n-k+1)}(t)\bigr).
\end{equation}
See~\cite[p.~139, Theorem~C]{Comtet-Combinatorics-74}. The Bell polynomials of the second kind $\bell_{n,k}$ satisfy
\begin{gather}\label{113-final-formula}
\sum_{n=k}^\infty \bell_{n,k}(x_1,x_2,\dotsc,x_{n-k+1})\frac{t^n}{n!}
=\frac1{k!}\Biggl(\sum_{m=1}^\infty x_m\frac{t^m}{m!}\Biggr)^k,\\
\label{Bell(n-k)}
\bell_{n,k}\bigl(abx_1,ab^2x_2,\dotsc,ab^{n-k+1}x_{n-k+1}\bigr)
=a^kb^n\bell_{n,k}(x_1,x_n,\dotsc,x_{n-k+1}),\\
\label{Bell-mions-recurs}
\bell_{n,k}\biggl(\frac{x_2}2, \frac{x_3}3,\dotsc,\frac{x_{n-k+2}}{n-k+2}\biggr)
=\frac{n!}{(n+k)!}\bell_{n+k,k}(0,x_2,\dotsc,x_{n+1}),
\end{gather}
where $a$ and $b$ are any complex numbers. See~\cite[pp.~133 and~135--136]{Comtet-Combinatorics-74}.
\par
Let
\begin{equation}\label{g(x)-dfn}
h(x)=\int_0^\infty\biggl(\int_{1/e}^1 t^{xu-1}\td t\biggr)e^{-u}\td u.
\end{equation}
It is clear that, for $\ell\in\mathbb{N}$,
\begin{multline*}
h^{(\ell)}(x)=\int_0^\infty\biggl[\int_{1/e}^1 t^{xu-1}(\ln t)^\ell\td t\biggr]u^\ell e^{-u}\td u\\
\to\int_0^\infty\biggl[\int_{1/e}^1 \frac{(\ln t)^\ell}t\td t\biggr]u^\ell e^{-u}\td u
=\frac{(-1)^\ell\ell!}{\ell+1}
\end{multline*}
as $x\to0$.
Applying in~\eqref{Bruno-Bell-Polynomial} $f(v)=v^k$ and the function~\eqref{g(x)-dfn} to compute~\eqref{n-times-diriv} reveals
\begin{multline}\label{s(n-k)-dinal-eq}
s(n,k)=\binom{n}{k}\lim_{x\to0}\sum_{m=1}^{n-k}f^{(m)}(h(x)) \bell_{n-k,m}\bigl(h'(x),\dotsc,h^{(n-k-m+1)}(x)\bigr)\\
=\begin{cases}\displaystyle
\binom{n}{k}\lim_{x\to0}\sum_{m=1}^{k}f^{(m)}(h(x)) \bell_{n-k,m}\bigl(h'(x),\dotsc,h^{(n-k-m+1)}(x)\bigr), & n>2k\\ \displaystyle
\binom{n}{k}\lim_{x\to0}\sum_{m=1}^{n-k}f^{(m)}(h(x)) \bell_{n-k,m}\bigl(h'(x),\dotsc,h^{(n-k-m+1)}(x)\bigr), & k\le n\le2k
\end{cases}\\
=\begin{cases}\displaystyle
\binom{n}{k}\sum_{m=1}^{k}f^{(m)}(h(0)) \bell_{n-k,m}\bigl(h'(0),\dotsc,h^{(n-k-m+1)}(0)\bigr), & n>2k\\ \displaystyle
\binom{n}{k}\sum_{m=1}^{n-k}f^{(m)}(h(0)) \bell_{n-k,m}\bigl(h'(0),\dotsc,h^{(n-k-m+1)}(0)\bigr), & k\le n\le2k
\end{cases}\\
=\begin{cases}
\begin{aligned}
\binom{n}{k}\sum_{m=1}^{k}&\frac{k!}{(k-m)!} \bell_{n-k,m}\biggl(-\frac{1!}2,\dotsc, \\
&\frac{(-1)^{n-k-m+1}(n-k-m+1)!}{n-k-m+2}\biggr), \quad n>2k;
\end{aligned}\\
\begin{aligned}
\binom{n}{k}\sum_{m=1}^{n-k}&\frac{k!}{(k-m)!} \bell_{n-k,m}\biggl(-\frac{1!}2,\dotsc, \\ &\frac{(-1)^{n-k-m+1}(n-k-m+1)!}{n-k-m+2}\biggr), \quad k\le n\le2k.
\end{aligned}
\end{cases}
\end{multline}
\par
Taking $x_m=\frac{m!}{m+1}$ in~\eqref{113-final-formula} and using~\eqref{gen-funct-3} give
\begin{gather*}
\sum_{n=k}^\infty \bell_{n,k}\biggl(\frac{1!}2,\frac{2!}3,\dotsc, \frac{(n-k+1)!}{n-k+2}\biggr)\frac{t^n}{n!}
=\frac1{k!}\Biggl(\sum_{m=1}^\infty \frac{t^m}{m+1}\Biggr)^k \\
=\frac{(-1)^k}{k!}\biggl[\frac{\ln(1-t)}{t}+1\biggr]^k
=\frac{(-1)^k}{k!}\sum_{i=0}^k\binom{k}i\biggl[\frac{\ln(1-t)}{t}\biggr]^i\\
=\frac{(-1)^k}{k!}\sum_{i=0}^k\binom{k}{i}\frac{i!}{t^i}\sum_{\ell=i}^\infty (-1)^\ell s(\ell,i)\frac{t^\ell}{\ell!}
=(-1)^k\sum_{i=0}^k\frac1{(k-i)!} \sum_{\ell=i}^\infty(-1)^\ell s(\ell,i)\frac{t^{\ell-i}}{\ell!}.
\end{gather*}
This implies that
\begin{align*}
\bell_{n,k}\biggl(\frac{1!}2,\frac{2!}3,\dotsc, \frac{(n-k+1)!}{n-k+2}\biggr)
&=n!(-1)^k\sum_{i=0}^k\frac{(-1)^{n+i}}{(k-i)!} \frac{s(n+i,i)}{(n+i)!}\\
&=(-1)^{n-k}\frac1{k!}\sum_{i=0}^k\frac{\binom{k}{i}}{\binom{n+i}{i}}(-1)^is(n+i,i).
\end{align*}
The formula~\eqref{Bell-Stir1st=eq} follows.
\par
Substituting~\eqref{Bell-Stir1st=eq} into~\eqref{Bell-mions-recurs} leads to
\begin{equation*}
\bell_{n+k,k}(0,1!,2!\dotsc,n!)=(-1)^{n-k}\binom{n+k}{k} \sum_{i=0}^k(-1)^i\frac{\binom{k}{i}}{\binom{n+i}{i}}s(n+i,i),
\end{equation*}
which may be rearranged as~\eqref{Bell-Stir1n!=eq}.
\par
By virtue of~\eqref{Bell(n-k)}, we have
\begin{multline}\label{Bell-simplify-eqminus}
\bell_{n-k,m}\biggl(-\frac12,\frac23,\dotsc, \frac{(-1)^{n-k-m+1}(n-k-m+1)!}{n-k-m+2}\biggr) \\
=(-1)^{n-k}\bell_{n-k,m}\biggl(\frac{1!}2,\frac{2!}3,\dotsc, \frac{(n-k-m+1)!}{n-k-m+2}\biggr).
\end{multline}
Substituting~\eqref{Bell-Stir1st=eq} into~\eqref{Bell-simplify-eqminus}, and then into~\eqref{s(n-k)-dinal-eq}, and simplifying find that
\begin{enumerate}
\item
when $2k\ge n\ge k\ge1$, we have
\begin{equation}\label{s(n-k)=s(n-k)-id-1}
s(n,k)=\sum_{m=1}^{n-k}\sum_{\ell=1}^m(-1)^{m+\ell} \binom{n}{k-\ell}\binom{k-\ell}{m-\ell} s(n-k+\ell,\ell);
\end{equation}
\item
when $n>2k>0$, we have
\begin{equation}\label{s(n-k)=s(n-k)-id-2}
s(n,k)=\sum_{m=1}^{k}\sum_{\ell=1}^m(-1)^{m+\ell} \binom{n}{k-\ell}\binom{k-\ell}{m-\ell} s(n-k+\ell,\ell).
\end{equation}
\end{enumerate}
Considering the convention that $s(n,k)=0$ for $0\le n<k$, we can unify the above two formulas~\eqref{s(n-k)=s(n-k)-id-1} and~\eqref{s(n-k)=s(n-k)-id-2} into
\begin{equation}
s(n,k)=\sum_{m=1}^{n}\sum_{\ell=1}^m(-1)^{m+\ell} \binom{n}{k-\ell}\binom{k-\ell}{k-m} s(n-k+\ell,\ell),
\end{equation}
which can be further formulated as~\eqref{s(n-k)=s(n-k)-id}.
\par
Interchanging two sums in~\eqref{s(n-k)=s(n-k)-id} and computing the inner sum yield
\begin{align*}
s(n,k)&=(-1)^k\sum_{\ell=k-n}^{k-1}(-1)^\ell\binom{n}{\ell}\Biggl[\sum_{m=k-\ell}^n(-1)^{m} \binom{\ell}{k-m}\Biggr] s(n-\ell,k-\ell)\\
&=(-1)^{n-k}\sum_{\ell=k-n}^{k-1}(-1)^\ell\binom{n}{\ell} \binom{\ell-1}{k-n-1}s(n-\ell,k-\ell)
\end{align*}
which may be rearranged as~\eqref{1stirling-diagonal-eq}.
The proof of Theorem~\ref{s(n-k)=Bell(n-k)-thm} is complete.

\section{Remarks}

\begin{rem}
The recurrence relations~\eqref{s(n-k)=s(n-k)-id} and~\eqref{1stirling-diagonal-eq} are neither ``triangular'', nor ``vertical'', nor ``horizontal'' recurrence relations as listed in~\cite[pp.~214--215, Theorems~A,~B, and~C]{Comtet-Combinatorics-74}, so we call them ``diagonal'' recurrence relations for the Stirling numbers of the first kind $s(n,k)$.
\end{rem}

\begin{rem}
The formula~\eqref{s(n-k)=s(n-k)-id} is also true if changing the sum over $m$ from $1$ to $k$ instead of from $1$ to $n$.
\end{rem}

\begin{rem}
Corollary~2.3 in~\cite{Filomat-36-73-1.tex} states that the Stirling numbers of the first kind $s(n,k)$ for $2\le k\le n$ may be computed by
\begin{equation}\label{s(n,k)-sum}
s(n,k)=(-1)^{n-k}(n-1)!\sum_{\ell_1=k-1}^{n-1} \frac1{\ell_1}\sum_{\ell_2=k-2}^{\ell_1-1}\frac1{\ell_2}\dotsm \sum_{\ell_{k-2}=2}^{\ell_{k-3}-1} \frac1{\ell_{k-2}} \sum_{\ell_{k-1}=1}^{\ell_{k-2}-1}\frac1{\ell_{k-1}}.
\end{equation}
This formula may be reformulated as
\begin{equation}
(-1)^{n-k}\frac{s(n,k)}{(n-1)!}= \sum _{m=k-1}^{n-1}\frac1m\biggl[(-1)^{m-(k-1)}\frac{s(m,k-1)}{(m-1)!}\biggr].
\end{equation}
\end{rem}

\begin{rem}
By applying the integral representation~\eqref{n-times-diriv}, some properties for the Stirling numbers of the first kind $s(n,k)$, including the logarithmic convexity with respect to $n\ge0$ of the sequence $\Bigl\{\frac{|s(n+k,k)|}{\binom{n+k}{k}}\Bigr\}_{n\ge0}$ for any fixed $k\in\mathbb{N}$, see~\cite[Corollary~5.1]{1st-Sirling-Number-2012.tex}, were established in~\cite[Section~5]{1st-Sirling-Number-2012.tex}.
\end{rem}

\begin{rem}
It is well known in combinatorics that
\begin{equation}\label{Bell-factorial-eq}
\bell_{n,k}(1!,2!,\dotsc,(n-k+1)!)=\binom{n}{k}\binom{n-1}{k-1}(n-k)!
\end{equation}
for $n\ge k\ge1$. See~\cite[p.~135, Theorem~B]{Comtet-Combinatorics-74}. We now recover this identity alternatively.
\par
In~\cite[Theorems~2.1 and~2.2]{exp-reciprocal-cm-IJOPCM.tex}, it was inductively obtained that, for $i\in\mathbb{N}$ and $t\ne0$,
\begin{equation}\label{exp-frac1x-expans}
\frac{\td^ie^{1/t}}{\td t^i}=(-1)^ie^{1/t}\frac1{t^{2i}}\sum_{k=0}^{i-1}a_{i,k}t^{k}
\end{equation}
and
\begin{equation}\label{g(t)-derivative}
\frac{\td^ie^{-1/t}}{\td t^i}=\frac{e^{-1/t}}{t^{2i}} \sum_{k=0}^{i-1}(-1)^ka_{i,k}{t^{k}},
\end{equation}
where
\begin{equation}\label{a-i-k-dfn}
a_{i,k}=\binom{i}{k}\binom{i-1}{k}{k!}
\end{equation}
for all $0\le k\le i-1$ and $a_{n,n-k}$ are Lah numbers $L(n,k)$. See also~\cite[Equations~(1.3) and~(1.4)]{QiBerg.tex}. For more information on Lah numbers $L(n,k)$, please refer to the recent references~\cite{DMST-MM-2013-Exp} and~\cite{Lindsay-Mansour-Shattuck-JComb-2011} and related reference therein.
\par
By~\eqref{Bruno-Bell-Polynomial} and~\eqref{Bell(n-k)}, it follows that, for $i\in\mathbb{N}$ and $t\ne0$,
\begin{equation}\label{exp-bruno-1}
\begin{split}
\frac{\td^ie^{1/t}}{\td t^i}&=e^{1/t}\sum_{k=1}^i \bell_{i,k}\biggl(-\frac{1!}{t^2},\frac{2!}{t^3},\dotsc,(-1)^{i-k+1}\frac{(i-k+1)!}{t^{i-k+2}}\biggr)\\
&=(-1)^ie^{1/t}\sum_{k=1}^i \frac1{t^{i+k}}\bell_{i,k}(1!,2!,\dotsc,(i-k+1)!)
\end{split}
\end{equation}
and
\begin{align*}
\frac{\td^ie^{-1/t}}{\td t^i}&=e^{-1/t}\sum_{k=1}^i \bell_{i,k}\biggl(\frac{1!}{t^2},-\frac{2!}{t^3},\dotsc,(-1)^{i-k}\frac{(i-k+1)!}{t^{i-k+2}}\biggr)\\
&=e^{-1/t}\sum_{k=1}^i (-1)^k \bell_{i,k}\biggl(-\frac{1!}{t^2},\frac{2!}{t^3},\dotsc,(-1)^{i-k+1}\frac{(i-k+1)!}{t^{i-k+2}}\biggr)\\
&=e^{-1/t}\sum_{k=1}^i \frac{(-1)^{i+k}}{t^{i+k}}\bell_{i,k}(1!,2!,\dotsc,(i-k+1)!).
\end{align*}
\par
Combining the formula~\eqref{exp-frac1x-expans} with~\eqref{exp-bruno-1} and the formula~\eqref{g(t)-derivative} with the above equation respectively show
\begin{equation*}
(-1)^i\frac1{t^{2i}}\sum_{k=0}^{i-1}a_{i,k}t^{k} =(-1)^i\sum_{k=1}^i \frac1{t^{i+k}}\bell_{i,k}(1!,2!,\dotsc,(i-k+1)!)
\end{equation*}
and
\begin{equation*}
\frac1{t^{2i}} \sum_{k=0}^{i-1}(-1)^ka_{i,k}{t^{k}}=\sum_{k=1}^i \frac{(-1)^{i+k}}{t^{i+k}}\bell_{i,k}(1!,2!,\dotsc,(i-k+1)!).
\end{equation*}
As a result,
\begin{equation*}
\sum_{k=1}^{i}a_{i,i-k}t^{k} =\sum_{k=1}^i \bell_{i,k}(1!,2!,\dotsc,(i-k+1)!) t^{k},
\end{equation*}
which implies
\begin{multline}
\bell_{n,k}(1!,2!,\dotsc,(n-k+1)!)=a_{n,n-k}\\*
=\binom{n}{n-k}\binom{n-1}{n-k}{(n-k)!} =\binom{n}{k}\binom{n-1}{k-1}(n-k)!,
\end{multline}
a recovery of the identity~\eqref{Bell-factorial-eq}.
\end{rem}

\begin{rem}
In~\cite{Guo-Qi-JANT-Bernoulli.tex} and~\cite{Special-Bell2Euler.tex} and related references therein, several special values of the Bell polynomials of the second kind $\bell_{n,k}$ are collected and applied.
\end{rem}

\begin{rem}
The term $(-1)^{\ell-1}$ in~\eqref{vert-rel-2} was misprinted as $(-1)^{n-1}$ in~\cite[p.~215, Theorem~B]{Comtet-Combinatorics-74}.
\end{rem}

\begin{rem}
This paper is a revised version of the preprint~\cite{notes-Stirl-No-JNT-rev.tex}.
\end{rem}

\end{document}